\documentclass[12pt]{article}
\usepackage{amsthm}
\usepackage[all]{xy}
\usepackage{enumerate}
\usepackage{amsfonts}
\usepackage{amsmath}

\def\ni{\noindent}
\def\sk1{\vskip 10pt}
\def\ds{\displaystyle}
\def\x{\times}

\newcommand{\R}{\mathbb{R}}
\newcommand{\bS}{\mathbb{S}}

\newcommand{\Z}{\mathbb{Z}}

\newcommand{\sbs}{\subset}
\newcommand{\ra}{\rightarrow}
\newcommand{\M}{{\cal{M}}}
\newcommand{\bM}{\bar{{\cal{M}}}}
\newcommand{\hM}{\hat{{\cal{M}}}}
\newcommand{\mP}{{\mathcal{P}}}
\newcommand{\mR}{{\mathcal{R}}}
\newcommand{\tM}{\tilde{{\cal{M}}}}
\newcommand{\tQ}{\tilde{Q}}
\newcommand{\tR}{\tilde{R}}
\newcommand{\p}{\partial}
\newcommand{\n}{\nabla}
\newtheorem{thm}{Theorem}
\newtheorem{lem}{Lemma}[section]
\newtheorem*{add}{Addendum}
\newtheorem*{prop}{Proposition}
\newtheorem*{Claim}{Claim}
\newtheorem{Cla}{Claim}

\theoremstyle{definition}
\newtheorem{rem}{Remark}
\newtheorem*{cav}{Caveat}
\newtheorem*{Rem}{Remark}
\newtheorem*{claim}{Claim}

\title{Harmonic Cellular Maps which are not Diffeomorphisms}
\author{F.T. Farrell and P. Ontaneda\thanks{The first author was
partially supported by a NSF grant.
The second author was supported in part 
by a research  grant from CNPq, Brazil.}}

\date{ } 

\begin{document}

\maketitle

\baselineskip = 20pt
\setcounter{section}{-1}
\section{Introduction}\label{S:sec0}

The use of harmonic maps has been spectacularly successful in proving rigidity (and
superrigidity) results for non-positively curved Riemannian manifolds.  This is witnessed for
example by results of Sui \cite{Sui}, Sampson \cite{Sa}, Corlette \cite{Co}, Gromov and
Schoen \cite{GS}, Jost and Yau \cite{JY}, and Mok, Sui and Yeung \cite{MSY}.  All of which are
based on the pioneering existence theorem of Eells and Sampson \cite{ES} and the uniqueness
theorem of Hartman \cite{Har} and Al'ber \cite{A}.  In light of this we believe that it is
interesting to demarcate this technique.  For example, it was shown in \cite{FJ2} that a
harmonic homotopy equivalence between closed negatively
curved Riemannian manifold is sometimes not a diffeomorphism,
even when one of the manifolds has constant sectional curvature equal to $-1$ (i.e.\ is a real
hyperbolic manifold).  Later other examples were given in \cite{FJO} and \cite{FOR} where such
a harmonic homotopy equivalence $f$ is not even a homeomorphism; even though the ones
constructed in \cite{FOR} are homotopic to diffeomorphisms.  In this paper we construct a
harmonic map $h$ between closed negatively curved Riemannian manifolds $M$ and $N$ which is
not a diffeomorphism but is the limit of a 1-parameter family of diffeomorphisms; 
in particular, $h$ is a cellular map.  In
our example either $M$ or $N$ (but of course not both) can be a real hyperbolic manifold and
the other have its sectional curvature pinched within $\epsilon$ of $-1$, where $\epsilon$ is
any preassigned positive number.  (We do not know whether such a harmonic map $h$ can ever be
a homeomorphism.  See our acknowledgment below.)  This result is contained in Theorem
\ref{T:thm1}, its Addendum and Theorem \ref{T:thm2}.  We construct such examples in all
dimensions $> 10$ and conjecture that this can be improved to all dimensions $\geq 6$.  

We have also discovered a curious relationship between the Poincar{\'e} Conjecture in low
dimensional topology and the existence of a certain type of harmonic map $k :M\to N$ between
high dimensional (i.e.\ dim $M >10$) closed negatively curved Riemannian manifolds.  Recall
that the Poincar{\'e} Conjecture asserts that the only simply connected closed 3-dimensional
manifold is the 3-sphere.  If this is true, then there exists such a harmonic map $k$ which is
homotopic to a diffeomorphism but cannot be approximated by homeomorphisms; i.e.\ is not a
cellular map.  (See \cite{E} for a discussion of cellular maps which are called cell like
maps in that article.)  In particular if $f : M\to N$ is any smooth map homotopic to $k$ and
$f_t$, $t\geq 0$, denotes the heat flow from $f_0=f$ to $f_\infty = k$ given by the
Eells-Sampson Theorem \cite{ES}, then $f_t$ is neither univalent (i.e.\ not one-to-one) nor an
immersion for all $t$ sufficiently large (i.e.\ all $t\geq T_f$ for some $T_f\in {\mathbb
R})$.
Again either $M$ or $N$ (but not both) can, in this example, be real hyperbolic and the
sectional curvatures of the other be pinched within $\epsilon$ of $-1$.  This result is
contained in the Addendum to Theorem \ref{T:thm2}.  

The key to these two theorems and their addenda is a Proposition.  We now describe this
Proposition and outline its proof.  Crucial use is made of the main result of \cite{FOR}.  In
that paper a pair of homeomorphic but not PL homeomorphic closed negatively curved Riemannian
manifolds $M$ and $\M$ are constructed satisfying:
\begin{enumerate}
\item  $M$ is real hyperbolic.
\item  $\M$ has a 2-sheeted cover $q : \hM \to \M$ where
$\hM$ admits a real hyperbolic metric $\nu$.
\end{enumerate}
Let $\mu$ be a given negatively curved Riemannian metric on $\M$ and $q^*(\mu)$ be the induced
Riemannian metric on $\hM$.  We would like to find a 1-parameter family of negatively curved
Riemannian metrics connecting $q^*(\mu)$ to $\nu$.  But we don't know how to do this.  In fact
this is in general an open problem; cf.\ \cite[Question 7.1]{BK}.  However by passing to a
large finite sheeted cover $r : \bM \to \hM$, we are able to connect $(q\circ r)^*(\mu)$ to
the real hyperbolic metric $r^*(\nu)$ by a 1-parameter family of negatively curved Riemannian
metrics; this is essentially the content of our Proposition in which $p=q\circ r$.  To
accomplish this, several results about smooth pseudo-isotopies are used; in particular, the
main result of \cite{FJ1} concerning the space of stable topological pseudo-isotopies of real
hyperbolic manifolds together with the comparison between the spaces of stable smooth and
stable topological pseudo-isotopies contained in \cite{BL} and \cite{Hat}.  And finally we need
Igusa's fundamental result \cite{I} comparing the spaces of pseudo-isotopies and stable
pseudo-isotopies.  We need that dim $M > 10$ in order to invoke Igusa's result.  We also need
a formula calculating the sectional curvatures of doubly warped products.  Although such a
formula is probably known to experts, we sketch a proof of it in the Appendix to this paper
for the sake of completeness.  

Finally our two theorems and their addenda are derived from our Proposition by using the
continuous dependence (in the $C^\infty$-topology) of the harmonic map homotopic to a homotopy
equivalence $f : (M,\mu_M)\to (N,\mu_N)$ on the negatively curved Riemannian metrics $\mu_M$
and $\mu_N$.  This dependence was proved by Sampson \cite{Sa}, Schoen and Yau \cite{SY}, and
Eells and Lemaire \cite{EL2}.  

In addition the derivation of the Addenda to Theorem \ref{T:thm2} depends on Scharlemann's
result \cite{Sc} which is the key unlocking the connection to the Poincar{\'e} Conjecture.
\sk1

\ni {\textbf{Acknowledgment.}}  We wish to thank David Gabai for asking whether a harmonic
homeomorphism between negatively curved Riemannian manifolds must be a diffeomorphism.

\section{Main results}\label{S:sec1}

In this section we state the main results of the paper which are Theorem \ref{T:thm1} and
\ref{T:thm2}, their Addenda, and the Proposition.  Then we show how the Proposition implies
the other main results.

\begin{thm}\label{T:thm1} For every integer $m>10$, there is a 
harmonic cellular map  $h:M_{1}\ra M_{2}$, between a pair
of closed negatively curved $m$-dimensional Riemannian
manifolds, which is not a diffeomorphism.
\end{thm}

\begin{add} The map $h$ in Theorem \ref{T:thm1} can be approximated 
by diffeomorphisms. Also,
 either $M_{1}$ or
$M_{2}$ can be chosen to be a real hyperbolic manifold and
the other chosen to have its sectional curvatures pinched within 
$\epsilon$ of -1; where $\epsilon$ is any preassigned positive
number.
\end{add}

\begin{rem}\label{R:rem1}  Siebenmann \cite{Si} showed that a continuous map $f : X\to Y$
between a pair of closed manifolds of dimension $\geq 5$ is cellular if and only if it is the
limit of homeomorphisms.
\end{rem}

\begin{rem}\label{R:rem2}  Note that a smooth non-diffeomorphic cellular map $h$ between
closed smooth manifolds cannot be a smooth immersion; i.e.\ $dh$ is not one-to-one on some
tangent space.  However we do not know whether the map $h$ we construct in proving Theorem
\ref{T:thm1} is univalent (i.e.\ one-to-one) or whether a harmonic homeomorphism between
closed negatively curved Riemannian manifolds must always be a diffeomorphism.  
\end{rem}

\begin{thm}\label{T:thm2} For every integer $m>10$, and
$\epsilon >0$, there are
an $m$-dimensional closed orientable smooth
manifold $\M$, and a
$C^{\infty}$ family of Riemannian metrics $\mu_{s}$, 
on $\M$, $s\in [0,1]$, such that:
\begin{enumerate}[{\rm(a)}]
\item $\mu_{1}$ is hyperbolic. 
\item The sectional curvatures of $\mu_{s}$, $s\in [0,1]$, 
are all in interval $(-1-\epsilon, -1+\epsilon)$.
\item  The maps $k$ and $l$ are both not univalent (i.e.\ not one-to-one) where $k :
(\M,\mu_0) \to (\M,\mu_1)$ and $l : (\M,\mu_1) \to (\M,\mu_0)$ are the unique harmonic maps
homotopic to \text{\rm id}$_\M$.
\end{enumerate}
\end{thm}

\begin{add}  Assuming that the Poincar{\'e} Conjecture is true, then the harmonic maps $k$ and
$l$ (of Theorem \ref{T:thm2}) are not cellular.  And consequently the maps $k_t$ and $l_t$ in
the heat flow of\break {\rm id} $= k_0$ to $k = k_\infty$ and of {\rm id} $= l_0$ to $l =
l_\infty$ are not univalent for all $t$ sufficiently large.
\end{add}

\begin{rem}\label{R:rem3}  The heat flow $k_t$ mentioned in this Addendum refers to the
solution to the initial value problem.
\begin{equation*}
\frac{\partial k_t}{\partial t} = \tau(k_t),\ \ \ k_t\big|_{t=0} = \text{id}
\end{equation*}
where $\tau(k_t)$ is the tension field of $k_t$.  A fundamental result due to Eells and Sampson
\cite{ES} is that this PDE has a unique solution $k_t$ (for all $t\geq 0$) and that
$\ds{\lim_{t\to\infty}}k_t = k$; cf.\ \cite[pp.\ 22-24]{EL4}.  
\end{rem}

These theorems and their addenda are a consequence of the following result.

\begin{prop}
Given an integer $m>10$ and a positive number $\epsilon$, there
exist a $m$-dimensional closed orientable real hyperbolic manifold
$M$ and a smooth manifold $\M$ with the following properties:
\begin{enumerate}[{\rm(i)}]
\item $M$ is  homeomorphic to $\M$.
\item $M$ is not PL homeomorphic to $\M$.
\item $\M$ admits a Riemannian metric $\mu$, whose sectional 
curvatures are all in the interval $(-1-\epsilon , -1+\epsilon )$.
\item There is a finite sheeted cover  $p:\bM\ra\M$ and a 
one-parameter $C^{\infty}$ family of Riemannian metrics $\mu_{s}$, 
on $\bM$, $s\in [0,1]$, such that $\mu_{0}=p^{*}\mu$ and $\mu_{1}$ is
hyperbolic. The sectional curvatures of $\mu_{s}$, $s\in [0,1]$, 
are all in the interval $(-1-\epsilon, -1+\epsilon)$.
\end{enumerate}
\end{prop}

Before deducing these two theorems and their addenda from the Proposition, we recall some
results concerning harmonic maps and introduce some notation.

Let $X$ and $Y$ be closed negatively curved 
Riemannian manifolds with Riemannian metrics 
$\mu_{X}$ and $\mu_{Y}$, respectively. Let $g:X\ra Y$ be a homotopy
equivalence. Then there is a unique harmonic map $k:X\ra Y$
homotopic to $g$, given by the fundamental existence result of Eells 
and Sampson
\cite{ES} and uniqueness result by Hartman \cite{Har} and Al'ber \cite{A}.
$k$ depends on $\mu_{X}$, $\mu_{Y}$ and $g$. We write
$k=har(\mu_{X} ,\mu_{Y},g)$. In fact, fixing $g$, the map
\begin{equation*}
\begin{aligned}
har_g : Met^{(-)}(X) \x Met^{(-)}(Y) &\to C^\infty(X,Y)\\
(\mu_X,\mu_Y) &\mapsto har(\mu_X,\mu_Y,g)
\end{aligned}
\end{equation*}
is continuous because of \cite{Sa}, \cite{SY} and \cite{EL2}; cf.\ \cite[\S2.18]{EL4}.  Here
$Met^{(-)}(\cdot)$ is the space of negatively curved smooth Riemannian metrics,
with the $C^{\infty}$ topology. Note that $Met^{(-)}(\cdot)$ is an
open set of the space $Met(\cdot)$
of smooth Riemannian metrics, with the $C^{\infty }$ topology.

\begin{proof}[Proof of Theorems \ref{T:thm1} and \ref{T:thm2}]
We prove the theorems assuming the Proposition.

Let $M$, $\M$, $\bM$, $p$, $\mu$ and $\mu_{s}$ be  as in the Proposition.
Let $f:M\ra\M$ be a homeomorphism. Write $\mu_{M}$ for the hyperbolic
metric of $M$.
Let $k:M\ra\M$ be the unique harmonic map homotopic to $f$, where we
consider $M$ with metric $\mu_{M}$ and $\M$ with metric $\mu$.
Using the notation above, we have $k=har(\mu_{M} ,\mu ,f)$.

Let $q:\bar{M}\to M$ be the pullback of $p:\bM\to\M$ via 
$f:M\to\M$ and $\bar f : \bar M \to \bM$ be the  
lifting of $f$ occurring in the pullback diagram 
\begin{equation*}
\xymatrix{\bar M \ar@{->}[r]^-{\bar f} \ar@{->}[d]_-q &\bM\ar@{->}[d]^-p\\
M \ar@{->}[r]^-f &\bM.
}
\end{equation*}
Write $\mu_{\bar{M}}$ for the hyperbolic metric $q^{*}(\mu_{M})$
on $\bar{M}$. On $\bM$ define
\begin{equation*}
h_{s} = har_{\bar{f}}(\mu_{\bar M},\mu_{s})
= har(\mu_{\bar M},\mu_{s}).
\end{equation*}

Then $s\mapsto h_{s}$ is a continuous map from [0,1] to 
$C^{\infty }(\bar{M},\bM )$.

Let $\bar{k}:\bar{M}\ra\bM$ be the lifting of $k:M\ra\M$. It is 
easily deduced from \cite[2.2.0 and 2.3.2]{EL3} that $\bar{k}$ 
is also a harmonic map
from $\bar{M}$, with metric $\mu_{\bar M}=q^{*}(\mu_{M})$, to   $\bM$, 
with metric $\mu_{0}=p^{*}(\mu)$. Note that $\bar{k}$ is homotopic to
$\bar{f}$.

But we also have  that $h_{0}$ is the unique harmonic map homotopic to
$\bar{f}:\bar{M}\ra\bM$, where we consider $\bar{M}$ with metric
$\mu_{\bar M}$ and $\bM$ with metric $\mu_{0}$.
Hence $h_{0}=\bar{k}$.

\begin{Claim} $h_{0}$ is not univalent.
\end{Claim}

We use the same argument as in \cite[pp.\ 229-230]{FOR}.

It is enough to prove that $k$ is not univalent.
Suppose $k$ is univalent. Then $k$ is a $C^{\infty}$-homeomorphism
between $M$ and $\M$, and hence $M$ and $\M$ are $PL$ homeomorphic,
by the $C^{\infty}$-Hauptvermutung proven by M. Scharlemann and
L. Siebenmann \cite{SS}. This contradicts (ii) of the Proposition
and proves the claim.

Note that $h_{1}$ is a harmonic map between the hyperbolic manifolds
$\bar M$ and $\bM$ (with metrics $\mu_{\bar M}$ and $\mu_{1}$, 
respectively) homotopic to the homotopy equivalence
$\bar{f}:\bar{M}\ra\bM$. Hence, by Mostow's Rigidity Theorem
\cite{M}, $h_{1}$ is an isometry. In particular $h_{1}$ is a
diffeomorphism.

We have proven that there is a continuous map $s\mapsto h_{s}$
from [0,1] to $C^{\infty }(\bar{M},\bM )$, with the following
properties:
\begin{enumerate}[(a)]
\item $h_{0}$ is not univalent.
\item $h_{1}$ is a diffeomorphism.
\item $h_{s}$ is a harmonic map between the hyperbolic
manifold $\bar{M}$ and the negatively curved Riemannian
manifold $\bM$ (with metric $\mu_{s}$).
\end{enumerate}

Define 
\begin{equation*}
s_{0}=\inf\{s\in[0,1] : h_{s}\ \text{is a diffeomorphism}\}.
\end{equation*} 

Because the space of diffeomorphisms, from $\bar{M}$ to $\bM$, is open 
in $C^{\infty }(\bar{M},\bM )$,  we have that $s_{0}<1$. Also,
since $h_{0}$ is not univalent, $h_{s_{0}}$ is {\it not} a diffeomorphism.
Moreover, $h_{s_{0}}$ is a cellular map, since it can be approximated
by the diffeomorphisms $h_{t}$, $t\in (s_{0},1]$.

Take $h=h_{s_{0}}$, $M_{1}=\bar{M}$ with the hyperbolic metric
$\mu_{\bar{M}}$, and  $M_{2}=\bM$ with the negatively curved metric
$\mu_{s_{0}}$. Then the harmonic map $h:M_{1}\ra M_{2}$ is a 
cellular map which is not a diffeomorphism. The sectional
curvatures of $M_{2}$ lie in the interval $(-1-\epsilon , -1+\epsilon )$.
This proves Theorem \ref{T:thm1} and part of its Addendum. To prove that we can
take $M_{2}$ to be hyperbolic, just repeat the argument above with
$f^{-1}:\M\ra M$, with the obvious modifications. This completes the
proof of Theorem \ref{T:thm1} and its Addendum.
\end{proof}

To prove Theorem 2, let the manifold $\M$ in it be the manifold $\bM$ of the Proposition and
let the Riemannian metrics $\mu_s$ in Theorem \ref{T:thm2} be those of the Proposition.  Then
notice that $l$ is $h_0\circ h^{-1}_l$; which is harmonic since $h_1$ is an isometry.  Also
$k$ is $\hat h^{-1}_1 \circ \hat h_0$ where 
\begin{equation*}
\hat h_t = har_{\bar f^{-1}}(\mu_t,\mu_{\bar M}).
\end{equation*}
And as mentioned above, it can be shown analogously that $\hat h_1$ is an isometry while $\hat
h_0$ is not univalent.

To prove the Addendum to Theorem 2, it suffices to show that $h_0$ and $\hat h_0$ are both not
cellular.  We will only explicitly show this for $h_0$ since the argument for $\hat h_0$ is
analogous.

\begin{cav} In this argument we now revert to our earlier notation given on line 5 in the
``Proof of Theorems \ref{T:thm1} and \ref{T:thm2}'' where 
\begin{equation*}
k = har(\mu_M,\mu,f);
\end{equation*}
$k$ will no longer denote $har(\mu_0,\mu_1,\text{id})$.  
\end{cav}

Recall that we are now assuming that the Poincar{\'e} Conjecture is true.  Under this
assumption Scharlemann's main result in \cite{Sc} and his result with Siebenmann \cite{SS}
routinely combine to the sharper statement that the smooth map $k$ is not cellular since $M$
and $\M$ are not PL-homeomorphic.  Now recall that a continuous map between closed manifolds
is cellular if and only if the inverse image, under it, of each contractible open subset of
the range is contractible; cf.\ \cite{E}.  Hence there exists a contractible open subset $U$
of $\M$ such that $k^{-1}(U)$ is not contractible.  Now consider the open subset $W$ of $\bar
M$ defined by
\begin{equation*}
W = q^{-1}(k^{-1}(U)) = \bar k^{-1}(p^{-1}(U)).
\end{equation*}
The diagram 
\begin{equation*}
\xymatrix{
W \ar@{->}[r]^-{\bar k} \ar@{->}[d]_-q &p^{-1}(U) \ar@{->}[d]^-p\\
k^{-1}(U) \ar@{->}[r]_-k &U
}
\end{equation*}
shows that $q : W\to k^{-1}(U)$ is the pullback of the covering space $p : p^{-1}(U) \to U$.
But $p : p^{-1}(U) \to U$ is a trivial covering space since $U$ is contractible; consequently,
$q : W \to k^{-1}(U)$ is also trivial.  Now let $V$ be a sheet of $p : p^{-1}(U) \to U$, then
$V$ is a contractible open subset of $\bM$ since it is homeomorphic to $U$.  And one easily
sees that $\bar k^{-1}(V)$ is homeomorphic to $k^{-1}(U)$; consequently, $h_0 = \bar k$ is not
cellular.  This completes the proof of Theorem \ref{T:thm2} and its Addendum.

We will prove the Proposition in section 3.  In the next section (section \ref{S:sec2}) we
give three lemmas which will be needed to prove the Proposition.

\section{Preliminary Lemmas}\label{S:sec2}

The first lemma we state is similar to the Lemma of \cite{FOR}.

\begin{lem}\label{L:lem21} Given an integer $m\geq 6$ and a 
positive number  $r$,
there exist closed connected oriented real hyperbolic manifolds 
$M$, $N$, $T$ and a pair of cohomology classes $\alpha \in 
H^{1}(M,\Z_{2})$ and $\beta \in H^{2}(M,\Z_{2})$
satisfying the following properties:
\begin{enumerate}[{\rm(1)}]
\item {\rm dim}$(M)=m$ and $T$ is a totally geodesic codimension-one
submanifold of  $M$.
\item $N$ is a totally geodesic framable codimension-two submanifold 
of $M$, whose
normal geodesic tubular neighborhood has width $\geq r$.
\item The isometry class of $N$ depends only on $m$ (not on $r$).
\item $\alpha \cup \beta \neq 0$.
\item $\alpha $ is the Poincare dual of the homology class represented
by $T$ in $H_{m-1}(M,\Z_{2})$.
\item  $\beta $ is the Poincare dual of the homology class represented
by $N$ in $H_{m-2}(M,\Z_{2})$.
\end{enumerate}
\end{lem}

\begin{proof}
Our proof is the same as the proof of the Lemma in \cite{FOR},
just interchange $n_{1}$ and $n_{2}$ at the beginning of that proof.
\end{proof}

We now give a geometric lemma, but first we introduce some notation 
and make some comments.

Let $M$ be a Riemannian manifold, with Riemannian metric $\sigma$.
Let also $\phi :\R\ra (0,\infty )$ be a smooth function.
Consider the warped metric $\rho = \phi^{2}\sigma+ dt^{2}$ on
$M\x \R$. A classic formula of Bishop and O'Neill \cite{BO},
gives the sectional curvatures $K_{\rho}$, of the Riemannian metric
$\rho,$ in terms of $\phi$ and the sectional curvatures 
$K_{\sigma}$ of $\sigma$:
\begin{equation*}
K_{\rho}(P)= -\frac{\phi''(t)}{\phi 
(t)} s^{2}+ \left( \frac{K_{\sigma}(u,v)-(\phi'(t))^{2}}{\phi^{2}(t)}\right)\| u 
\|^{2}
\end{equation*}

Here $P\subset T_{(x,t)}(M\x\R )= T_{x}M\oplus\R $ is the two-plane 
generated by the orthonormal basis $\{ u+s\frac{\partial}{\partial 
t},v\}$, where $u,v\in T_{x}M$. Note that  $s^{2}+\| u\|^{2}=1$,
$\| v\|^{2}=1$ and $\langle u,v\rangle =0$. It follows that $K_{\rho}(P) $ is a convex
linear combination of $- \frac{\phi''(t)}{\phi (t)}$ and 
$\frac{K_{\sigma}(u,v)-(\phi '(t))^{2}}{\phi^{2}(t)}$.

We consider now doubly warped metrics. Let $M_{1}$ and $M_{2}$
be Riemannian manifolds with Riemannian metrics $\sigma_{1}$ and 
$\sigma_{2}$, respectively. Let  also $\phi_{i} :\R\ra (0,\infty )$ 
be  smooth functions, $i=1,2.$ Define the doubly warped 
metric $\rho$ on $M_{1}\x M_{2}\x \R$ by
\begin{equation*}
\rho = \phi_{1}^{2}\sigma_{1}+ \phi_{2}^{2}\sigma_{2}+ dt^{2}
\end{equation*}

A generalization of Bishop-O'Neill's formula 
gives the sectional curvatures $K_{\rho}$ of the Riemannian metric
$\rho$ in terms of $\phi_{i}$ and the sectional curvatures 
$K_{i}$ of $\sigma_{i}$, $i=1,2$ (see the appendix for a proof):
\begin{equation*}
\begin{aligned}
K_{\rho}(P) &= -\frac{\phi_{1}''(t)}{\phi_{1} 
(t)}s^{2} \| u_{2}\|^2  -  \frac{\phi_{2}''(t)}{\phi_{2}(t)}s^{2}
\|  v_{2}\|^2\\ 
&+\ \left(\frac{K_{1}(u_{1},u_{2})-(\phi_{1}'(t))^{2}}{\phi_{1}^{2}(t)}\right)
(\| u_{1} \|^{2} \| u_{2} \|^{2} -\langle 
u_{1},u_{2}\rangle^{2})\\
&+\left(\frac{K_{2}(v_{1},v_{2})-(\phi_{2}'(t))^{2}}{\phi_{2}^{2}(t)}\right)
(\| v_{1} \|^{2} \| v_{2} \|^{2}  -\langle 
v_{1},v_{2}\rangle^{2})\\
&-\ \frac{\phi_{1}'(t)\phi_{2}'(t)}{\phi_{1}(t)\phi_{2}(t)}
( \| u_{1} \|^{2} \| v_{2} \|^{2}  +
\| v_{1} \|^{2} \| u_{2} \|^{2}
-2\langle u_{1},u_{2}\rangle \langle v_{1},v_{2}\rangle).
\end{aligned}
\end{equation*}

Here $P\subset T_{(x_{1}, x_{2},t)}(M_{1}\x M_{2}\x\R ) = 
T_{x_{1}}M_{1}\oplus T_{x_{2}}M_{2}\oplus \R $ is the two-plane 
generated by the orthonormal basis $\{ u_{1}+v_{1}+s\frac{\partial}{\partial 
t}, u_{2}+v_{2}\}$, where $u_{1},u_{2}\in T_{x_{1}}M_{1}$,
$v_{1},v_{2}\in T_{x_{2}}M_{2}$. 
Note that  $s^{2}+\| u_{1}\|^{2}+\| v_{1}\|^{2}=1$,
$\| u_{2}\|^{2}+\| v_{2}\|^{2}=1$, and $\langle u_{1},u_{2}\rangle 
+\langle v_{1},v_{2}\rangle =0$. It follows that $K_{\rho}(P) $ is a convex
linear combination of $- \frac{\phi_{i}''(t)}{\phi_{i} (t)}$,  
$\frac{K_{i}(\cdot,\cdot)-(\phi_{i}'(t))^{2}}{\phi_{i}^{2}(t)}$
and $- \frac{\phi_{1}'(t)\phi_{2}'(t)}{\phi_{1}(t)\phi_{2}(t)}$,
$i=1,2$.

Let $a,b\in\R$, $0<a<b$. Let also $\phi_{i} :\R\x (0,\infty )\to
(0,\infty) $, $i=1,2$ be two smooth functions. 
The first and second derivatives of $\phi_{i}$ with respect to the 
first variable, will be denoted by $\phi_{i}'$ and $\phi_{i}''$,
respectively.
In the next lemma we consider Riemannian 
metrics $\rho_{\alpha}$ on $M_{1}\x M_{2}\x [a,b]$:
\begin{equation*}
\rho_{\alpha}(x_{1},x_{2},t) = \phi_{1}^{2}(\alpha t,\alpha) \sigma_{1}(x_{1})+ 
\phi_{2}^{2}(\alpha t,\alpha ) \sigma_{2}(x_{2})+ \alpha^{2}dt^{2}
\end{equation*}

\begin{lem}\label{L:lem22} Let $M_{1}$ and $M_{2}$
be compact Riemannian manifolds with Riemannian metrics $\sigma_{1}$ and 
$\sigma_{2}$, respectively. Let $\phi_{i} :\R\x (0,\infty )\to 
(0,\infty) $, $i=1,2$ be two smooth functions. Suppose that
\begin{equation*}
\lim_{\alpha\to\infty}\frac{\phi_{i}'(\alpha t,\alpha )}{\phi_{i}(\alpha t,\alpha )}=
\lim_{\alpha\to\infty}\frac{\phi_{i}''(\alpha t,\alpha )}{\phi_{i}(\alpha 
t,\alpha )}=1
\end{equation*}
and 
\begin{equation*}
\lim_{\alpha\to\infty}\phi_{i}(\alpha t,\alpha )=\infty
\end{equation*}
$i=1,2$, uniformly for $t\in[a,b]$.
Then, given $\epsilon >0$, there is an $\alpha_{0}\in\R$, such that,
for all $\alpha >\alpha_{0}$, all sectional curvatures of $\rho_{\alpha}$
lie in the interval $(-1-\epsilon,-1+\epsilon)$.
\end{lem}

\begin{proof}
The manifold $M_{1}\x M_{2}\x [a,b]$ with the Riemannian metric $\rho_{\alpha}$ is 
isometric to $M_{1}\x M_{2}\x [\alpha a,\alpha b]$ with 
the doubly warped metric
\begin{equation*}
\bar{\rho}_{\alpha}(x_{1},x_{2},t) = \phi_{1}^{2}(t,\alpha 
)\sigma_{1}(x_{1})+ 
\phi_{2}^{2}( t,\alpha )\sigma_{2}(x_{2})+ dt^{2}.
\end{equation*}

(The isometry is $t\mapsto\alpha t$, $t\in [a,b]$.)

As mentioned above, the sectional curvatures of $\bar{\rho}_{\alpha}$
are  convex linear combinations of $-\frac{\phi_{i}''(t,\alpha )}
{\phi_{i} (t,\alpha )}$,  
$\frac{K_{i}(\cdot,\cdot)-(\phi_{i}'(t,\alpha ))^{2}}{\phi_{i}^{2}(t,\alpha )}$
and $-\frac{\phi_{1} '(t,\alpha )\phi_{2} '(t,\alpha )}{\phi_{1}(t,\alpha )
\phi_{2}(t,\alpha )}$,
$i=1,2$.

But $M_{1}$ and $M_{2}$ are compact, therefore the sectional curvatures
$K_{1}$ and $K_{2}$ are bounded. Consequently, choosing $\alpha_{0}$ 
sufficiently large, we can suppose that all the terms 
$-\frac{\phi_{i}''(t,\alpha )}
{\phi_{i} (t,\alpha )}$,  
$\frac{K_{i}(\cdot,\cdot)-(\phi_{i}'(t,\alpha ))^{2}}{\phi_{i}^{2}(t,\alpha )}$
and $-\frac{\phi_{1}'(t,\alpha )\phi_{2}'(t,\alpha )}{\phi_{1}(t,\alpha )
\phi_{2}(t,\alpha )}$,
$i=1,2$, are within $\epsilon$ of -1, for $\alpha > \alpha_{0}$ and 
$t\in [\alpha a,\alpha b]$. It follows that all
the sectional curvatures of $\bar{\rho}_{\alpha}$ (and 
$\rho_{\alpha}$) lie in the interval $(-1-\epsilon ,-1+\epsilon )$.
This proves the lemma.
\end{proof}

Finally, we will need the following lemma.

\begin{lem}\label{L:lem23} Let $N$ be a closed connected orientable real 
hyperbolic manifold of dimension $\geq 10$, and let 
$f:N\x [0,1]\ra N\x [0,1]$ be a diffeomorphism which is 
smoothly
pseudo-isotopic to the identity, rel $\partial$.
Then there is a positive integer $J$ and
a finite number of non-conjugate elements $a_{1},\ldots ,a_{l}\in \pi_{1}N$ with 
the following property:

If $\bar{p}:\bar{N}\ra N$ is a connected finite cover  
such that no conjugate of $a^J_i$ belongs to $\bar{p}_{*}(\pi_{1}\bar{N})\subset \pi_{1}N$, then
$\bar{f}:\bar{N}\x [0,1]\ra N\x [0,1]$ is smoothly
isotopic to the identity, rel $\partial$.
\end{lem}
 
\begin{Rem} Note that if $p':N'\ra N$ is a connected finite cover
that factors through $\bar{p}:\bar{N}\ra N$, then the lemma holds also 
for $p'$ and $N'$; that is, $f':N'\x [0,1]\ra N\x [0,1]$ is 
also smoothly isotopic to the identity.
\end{Rem}

\begin{proof}[Proof of Lemma \ref{L:lem23}] Let $n=\text{dim }N$, $P(\ )$ denote the space of
topological pseudo-isotopies and $\mP(\ )$ denote the space of stable topological
pseudo-isotopies.  Also let $P^{\text{diff}}(\ )$ denote the space of smooth pseudo-isotopies
and $\mP^{\text{diff}}(\ )$ the space of stable smooth pseudo-isotopies.  Recall that we have
canonical stabilization maps 
\begin{equation*}
\begin{aligned}
\iota^{\text{diff}} &: P^{\text{diff}}(\ ) \to \mP^{\text{diff}}(\ )\ \ \text{ and}\\
\iota &: P(\ ) \to \mP(\ )
\end{aligned}
\end{equation*}
such that the following square of maps commutes:
\begin{equation*}
\xymatrix{
P^{\text{diff}}(\ ) \ar@{->}[r]^-{\iota^{\text{diff}}} \ar@{->}[d] &\mP^{\text{diff}}(\ )
\ar@{->}[d]\\
P(\ ) \ar@{->}[r]^-\iota &\mP(\ ).
}
\end{equation*}
The vertical arrows in this square denote the natural forget structure maps.  It is a
consequence of work of Burghelea and Lashof \cite{BL} and Cerf \cite{Ce} that the forgetful
map $\mP^{\text{diff}}(\ ) \to \mP(\ )$ induces an isomorphism on $\pi_0$; cf.\ \cite[p.\
12]{Hat}.  Furthermore $\iota$ and $\iota^{\text{diff}}$ induce isomorphisms on $\pi_0$ for all
manifolds of dim $> 10$ by Igusa \cite{I}.  Consequently the square shows that the forgetful
map $P^{\text{diff}}(\ ) \to P(\ )$ also induces an isomorphism on $\pi_0$ for all smooth
manifolds of dim $> 10$. 

Now these four isomorphisms combined with \cite[Theorem 6.0]{FJ1} show that we can assume that
the pseudo-isotopy of $f$ (given in Lemma \ref{L:lem23}) is supported on the disjoint
neighborhoods of a finite number of (non-conjugate and non-trivial) embedded loops in $N$.
These neighborhoods are diffeomorphic to ${\mathbb D}^n \x \bS^1$, where ${\mathbb D}^n$ is the
closed $n$-disc.  Hence the lemma follows from the following claim.
\end{proof}

\begin{claim}  For every smooth pseudo-isotopy $F : [0,1] \x {\mathbb D}^n \x \bS^1 \to [0,1] \x
{\mathbb D}^n \x \bS^1$, there is a $J$ such that if $\bar F$ is the lifting of $F$ by the
connected $j$-sheeted cover $[0,1] \x {\mathbb D}^n \x \bS^1 \to [0,1] \x {\mathbb D}^n\x 
\bS^1$, with $j\geq J$,
then $\bar F$ is smoothly isotopic to the identity, rel$(\{0\}\x {\mathbb D}^n\x \bS^1 \cup
[0,1] \x \partial{\mathbb D}^n\x \bS^1)$.  
\end{claim}

\begin{proof}[Proof of Claim]  Because of the above discussion it suffices to show that $\bar
F$ is topologically isotopic to the identity, rel$(\{0\}\x {\mathbb D}^n\x \bS^1\cup [0,1]\x
\partial{\mathbb D}^n\x \bS^1)$.

Given $\epsilon > 0$, by taking $j$ sufficiently large, we have that $\bar F$ becomes
$\epsilon$-controlled in the $\bS^1$-direction.  Then by appropriately shrinking
inwards in the $[0,1]\x {\mathbb D}^n$ direction (an Alexander type isotopy) we get
$\epsilon$-control in all directions.  Hence, by \cite{FJ2}, $\bar F$ can be topologically
isotoped to the identity, rel$(\{0\}\x {\mathbb D}^n\x \bS^1 \cup [0,1] \x
\partial{\mathbb D}^n\x \bS^1)$.  

This proves the claim and thus completes the proof of Lemma \ref{L:lem23}.
\end{proof}
 
\section{Proof of the Proposition}\label{S:sec3}

Let $M$, $N$ and $T$ be as in Lemma \ref{L:lem21}, relative to a sufficiently 
large positive real number $r$. (How large is sufficient, will 
presently become clear.) Define $P=N\cap T$. 
Since $\alpha \cup\beta\neq 0$, $N$ and $T$ intersect transversally 
and  $P$ is a codimension-three totally 
geodesic submanifold of $M$. Moreover  $\alpha \cap\beta$
is the Poincare dual of the cycle represented by $P$ in
$H_{m-3}(M,\Z_{2})$.

Since  the trivial normal geodesic tubular neighborhood of $N$ has 
width $\geq r$, we can identify the tubular neighborhood $V$, of width $r$,
with $N\x B$, where $B\subset\R^{2}$ is the open ball, centered at 
the origin, of radius $r$. This identification is a metric 
identification on
\begin{equation*}
\begin{aligned}
V-N &= (N\x B) - (N\x \{0\})\\
&= N\x (B-\{0\})\\
&= N\x \bS^1\x (0,r),
\end{aligned}
\end{equation*}
where we consider $N\x \bS^1 \x (0,r)$ with the doubly warped Riemannian metric
\begin{equation*}
\rho(x,u,t)=\cosh^{2}(t) \sigma_{N}(x)+\sinh^{2}(t) \sigma_{\bS^{1}}(u)+dt^{2}
\end{equation*}
Here $\sigma_{N}$ is the hyperbolic metric on $N$ and $\sigma_{\bS^{1}}$
is the canonical Riemannian metric on 
\begin{equation*}
\bS^1 = \{(x,y) \in {\mathbb R}^2\mid x^2+y^2 = 1\}.  
\end{equation*}
Define $N_0$, $P_0$, $Q$ and $R$ by
\begin{equation*}
\begin{aligned}
N_0 &= N\x \{(1,0)\} \x \left\{\frac{r}{2}\right\}\\
&\subset N\x \bS^1\x (0,r) = (V-N)\\
&\subset M\\
P_0 &= N_0 \cap T\\
Q &= N\x \bS^1\x \left\{\frac{r}{2}\right\}\\
R &= N\x \bS^1_+ \x \left\{\frac{r}{2}\right\}\ \text{ where}\\
\bS^1_+ &= \{(x,y) \in S^1\mid x > 0\}.
\end{aligned}
\end{equation*}
It follows from these definitions that 
\begin{equation*}
\begin{aligned}
P_0 &= P\x \{(1,0)\}\x \left\{\frac{r}{2}\right\}\\
P_0 &\subset N_0 \subset Q\subset V\subset M\\
R &\text{ is diffeomorphic to }\ N_0\x [0,1].
\end{aligned}
\end{equation*}

The smooth manifold $\M$ of the statement of the Proposition 
is constructed by cutting $M$ apart along $Q$ and gluing back with a 
twist $f:Q\ra Q$. For the details of this construction see the ``Proof 
of Corollary'' in \cite[pp.\ 230-233]{FOR} with the following 
modifications:
\begin{itemize}
\item Replace $N$ in \cite{E} by $Q$.
\item Replace $N$ by $T$ in the definition of $\hat{\alpha}$ at the top of 
\cite[p.\ 232]{E}.
\end{itemize}

Then we have that this smooth manifold $\M$, constructed as above, satisfies
the following properties:
\begin{enumerate}
\item There is a homeomorphism $g:\M\to M$.
\item $\M$ is not $PL$ homeomorphic to $M$.
\item There is a two-sheeted connected double cover $q': \tilde{M}\to M$
such that the lifting $\tilde{g}:\tM\to\tilde{M}$, of $g$, is 
(topologically) pseudo-isotopic to a diffeomorphism.
\end{enumerate}
Here $p' : \tM \to \M$ is the double cover induced from $q'$ via $g$ which occurs in the
pullback diagram 
\begin{equation*}
\xymatrix{
\tM \ar@{->}[r]^-{\tilde g} \ar@{->}[d]_-{p'} &\tilde M \ar@{->}[d]^-{q'}\\
\M \ar@{->}[r]^-g &M.
}
\end{equation*}

Let $\tilde{N_{0}}$, $\tR$, $\tQ$, be the liftings, by $q'$, of 
$N_{0}$, $R$, $Q$, respectively.
It is clear from the construction of $q'$ that $\tQ$ is connected.
Note that $\tR$ is also diffeomorphic to $\tilde{N_{0}}\x [0,1]$.
Let $\tilde{f}:\tQ\ra\tQ$ be the lifting of $f:Q\ra Q$.

\begin{Claim} We can choose $f$ such that:
\begin{enumerate}[{\rm(1)}]
\item $f$ is the identity outside $R\sbs Q$.
\item $\tilde{f}$ is smoothly pseudo-isotopic to the identity.
\end{enumerate}
\end{Claim}

\begin{proof}[Proof of this Claim]
To prove (1) just note that the map $\gamma$ in \cite{E}
can be chosen to be constant outside a small neighborhood of
$P_{0}$ in $M$. Hence $f$ can chosen to be the identity outside
a small neighborhood of $P_{0}$ in $Q$. Since $P_{0}\sbs int (R)$,
it follows that $f$ can be chosen to satisfy (1) of the claim.

To prove (2) consider the following diagram
\begin{equation*}
\xymatrix{
\tQ\x I \ar@{->}[r]^-{\tilde\sigma} \ar@{->}[d] &\tilde M \ar@{->}[rr]^-\xi \ar@{->}[d]_-{q'}
&&\bS^2 \x \bS^1 \ar@{->}[r]^-\psi \ar@{->}[d]^-{1_{\bS^2}\x p} &\bS^3\ar@{->}[d]^-\phi\\
Q\x I \ar@{->}[r]^-\sigma &M \ar@{->}[rr]^-{(\hat \beta\x \hat\alpha)\circ \Delta} &&\bS^2\x \bS^1
\ar@{->}[r]^-\psi &\bS^3 \ar@{->}[r]^-{\bar\eta} &Top/0
}
\end{equation*}
Here $Q\x I$ is a tubular neighborhood of $Q$ and $\sigma 
:Q\x  I \ra M$ is the inclusion. The diagram above is the
diagram of \cite[p.\ 233]{FOR} except for the first column and that the second vertical arrow
is now denoted $q'$ instead of $q$.
Since $\bar{\eta}\varphi \psi \xi \tilde{\sigma}$ is null
homotopic, we have that the differentiable structure on 
$\tQ\x I$, rel $\partial$, induced by the inclusion in
$\tM$ is smoothly concordant, rel $\partial$, to the one
induced by the inclusion in $\tilde{M}$. It follows that
$\tilde{f}$ is smoothly pseudo-isotopic to the identity. This 
proves the Claim.
\end{proof}

Since the fundamental group of $\tR$ injects into the fundamental
group of $\tilde{M}$, and $\pi_{1}\tilde{M}$ is residually finite,
we can apply Lemma \ref{L:lem23}  to $\tilde{f}|_{\tR}:\tR\ra\tR$ to obtain
 a finite cover $\bar{q}:\bar{M}\to\tilde{M}$ with $\bar{f}:\bar{Q}\to\bar{Q}$
smoothly isotopic to the identity.
Write $q=q'\circ {\bar q}:\bar{M}\to M$ and let $p:\bM\to\M$ be the covering space induced by
$g$ from $q$ via the following pullback diagram:
\begin{equation*}
\xymatrix{
\bM \ar@{->}[r]^-{\bar g} \ar@{->}[d]_-p &\bar M\ar@{->}[d]^-q\\
\M \ar@{->}[r]^-g &M.
}
\end{equation*}
So far we have obtained the following:
\begin{itemize}
\item $M$ is a closed connected orientable real hyperbolic  
manifold of dimension $>10.$
\item $\M$ is obtained from $M$ by cutting along the hypersurface $Q$ and 
gluing back with the twist $f:Q\ra Q$.
\item  $\M$ is homeomorphic, but not $PL$ homeomorphic, to $M$.
\item  There is a connected finite sheeted cover $p:\bar{M}\ra M$
such that $\bar{f}:\bar{Q}\ra\bar{Q} $ is smoothly isotopic to the identity 
$1_{\bar{Q}}$. In particular, $\bM$ is diffeomorphic to $\bar M$.
\end{itemize}

Note that $\bar{f}$ and $p$ depend only on $f$. This ends the 
topological part of the proof of the Proposition.

Let $\epsilon > 0$. We now  show how to construct the 
Riemannian metrics $\mu$ and $\mu_{s}$, of the  statement of the 
Proposition, with sectional curvatures in the interval 
$(-1-\epsilon , -1 +\epsilon )$.

Recall that $Q$ has the normal geodesic tubular neighborhood
in $M$,
\begin{equation*}
V-N = N\x \bS^1 \x (0,r)
\end{equation*}
which is equipped with the doubly warped Riemannian metric
\begin{equation*}
\rho(x,u,t)=\cosh^{2}(t)\sigma_{N}(x)+\sinh^{2}(t)\sigma_{\bS^{1}}(u)+dt^{2}
\end{equation*}
where $\sigma_{N}$ is the hyperbolic metric on $N$ and $\sigma_{\bS^{1}}$
is the canonical Riemannian metric on $\bS^{1}$.
Note that $N\x \bS^1 \x (0,r)$ with the Riemannian metric $\rho$ is isometric to $N\x \bS^1 \x
(0,6)$ with the Riemannian metric
\begin{equation*}
\rho_{r}(x,u,t)=\cosh^{2}(\alpha t)\sigma_{N}(x)+\sinh^{2}(\alpha t)
\sigma_{\bS^{1}}(u)+\alpha^{2}dt^{2}
\end{equation*}
where $\alpha = r/6$.

Let
\begin{equation*}
\begin{aligned}
\delta_{1} &: \R\to [-1,1]\\
\delta_{2} &: \R\to [0,1]\\
\delta_{3} &: \R\to [-1,1]
\end{aligned}
\end{equation*}
be smooth functions such that:
\begin{equation*}
\begin{aligned}
\delta_{1}(t) &=\begin{cases} 
-1 &\text{$t\leq 2$}\\
1 &\text{$t\geq 3$}\end{cases}\\
\delta_{2}(t) &=\begin{cases} 0 &\text{$t\leq 3$}\\
1 &\text{$t\geq 4$}\end{cases}\\
\delta_{3}(t) &= \begin{cases} 1 &\text{$t\leq 4$}\\
-1 &\text{$t\geq 5$}\end{cases}
\end{aligned}
\end{equation*}
and all $\delta_{i}$ are constant near 1,2,3,4,5.

Notice that $Q$ has the (smooth) tubular neighborhood
$[N\x \bS^{1}\x (0,r)]_{f}$, in $\M$
obtained from $N\x \bS^{1}\x (0,r)$ by cutting
along $Q=N\x \bS^{1}\x \{\frac{r}{2}\}$ and gluing back with
$f:Q\to Q$. Also note that $[N\x \bS^1 \x (0,r)]_f$ is diffeomorphic to
$[N\x \bS^{1}\x (0,6)]_{f}$ which is 
obtained from $N\x \bS^{1}\x (0,6)$ by cutting
along $Q=N\x \bS^{1}\x \{ 3\}$ and gluing back with
$f:Q\ra Q$.

On $[N\x \bS^{1}\x (0,6)]_{f}$, consider the following 
Riemannian metric:
\begin{equation*}
\lambda_{r}(x,u,t)= \begin{cases}
\rho_{r}(x,u,t) &\text{$t\leq 2$, $5\leq t$}\\ 
 & \\
\cosh^{2}(\alpha t) \sigma_{N}(x) + (\frac{e^{\alpha t}+
\delta_{1}(t)e^{-\alpha 
t}}{2})^{2} \sigma_{\bS^{1}}(u) + \alpha^2 dt^2 &\text{$2\leq t \leq 3$}\\
& \\
\cosh^{2}(\alpha t)\Bigl\{(1-\delta_{2}(t)) 
f^{*}[\sigma_{N}(x)+\sigma_{\bS^{1}}(u)]\\ 
+\ \delta_{2}(t) [\sigma_{N}(x) +\sigma_{\bS^{1}}(u)]\Bigr\}  +\alpha^{2}dt^{2}
&\text{$3\leq t\leq 4$}\\
& \\
\cosh^{2}(\alpha t) \sigma_{N}(x)  + (\frac{e^{\alpha t}+
\delta_{3}(t)e^{-\alpha t}}{2})^{2} \sigma_{\bS^{1}}(u) +\alpha^2 dt^2 &\text{$4\leq t\leq
5$}\end{cases}
\end{equation*}
where $\alpha =r/6$.

It can be verified from Lemma \ref{L:lem22} that, taking $r$ large enough,
$\lambda_{r}(x,u,t)$ has sectional curvatures within $\epsilon$ of $-1$,
for $t\leq 3$ and $4\leq t$. Also, for $3\leq t\leq 4 $, by taking $r$
large enough, $\lambda_{r}(x,u,t)$ has sectional curvatures  within 
$\epsilon$ of $-1$. (See \cite[pp.\ 11-13]{O}.)

We now define the metric $\mu$ on $\M$ in the following way
(see \cite{O} for more details):
$\mu (p)$ is the hyperbolic metric, for $p\notin [N\x \bS^{1}\x 
(r/3,5r/6)]_{f}$, and $\mu (p)$ is the pullback of the metric 
$\lambda_{r}$ by the map $t\mapsto t/\alpha $,
where $\alpha =r/6, $ and $p=(x,u,t)\in  [N\x \bS^{1}\x 
(0,r)]_{f}\sbs\M$.

Let $s\mapsto f_{s}$, $s\in [0,1/2]$ be a smooth isotopy of $\bar{f}$, with
$f_{0}=\bar{f}$ and $f_{1/2}=1_{\bar{Q}}$. We assume that this isotopy
is constant near $0$ and $1/2$.
Let $\M_{s}$,
be the smooth manifold obtained from $\bM$ by cutting along $\bar{Q}$
and gluing back with $f_{s}$.
We construct a family of Riemannian metrics
$\mu'_{s}$, $s\in [0,1/2]$, on $\M_{s}$.
This construction is identical to that of $\mu$,
just repeat all the
definitions and arguments  above,  writing a ``bar''
above each symbol and replacing $f$ by $f_{s}$. For instance,
on $[N\x \bS^{1}\x (0,6)]_{f_{s}}$
the Riemannian metric 
$(\bar{\lambda}_{r})_{s}$ is given by the following formula: 
\begin{equation*}
(\lambda_{r})_{s}(x,u,t)= \begin{cases}
\bar{\rho}_{r}(x,u,t) &\text{$t\leq 2$, $5\leq t$}\\ 
& \\
\cosh^{2}(\alpha t) \sigma_{\bar{N}}(x) + (\frac{e^{\alpha t}+
\delta_{1}(t)e^{-\alpha t}}{2})^{2} \sigma_{\bS^{1}}(u) + \alpha^2 dt^2&\text{$2\leq t\leq 3$}\\
& \\
\cosh^{2}(\alpha t)\Bigl\{ (1-\delta_{2}(t)) 
f^{*}_{s}[\sigma_{\bar{N}}(x)+\sigma_{\bS^{1}}(u)] \\
+\ \delta_{2}(t) [\sigma_{\bar{N}}(x) +\sigma_{\bS^{1}}(u)]\Bigr\}  +\alpha^{2}dt^{2}
&\text{$3\leq t\leq 4$}\\
& \\
\cosh^{2}(\alpha t) \sigma_{\bar{N}}(x)  +  (\frac{e^{\alpha t}+
\delta_{3}(t)e^{-\alpha 
t}}{2})^{2} \sigma_{\bS^{1}}(u)+\alpha^2dt^2 &\text{$4\leq t\leq 5$}\end{cases}
\end{equation*}
where $\alpha =r/6$.

In this way we obtain Riemannian metrics $\mu_{s}'$ on $\M_{s}$.

Since $f$, $\bar{f}$ (hence $f_{s}$) do not depend on $r$,
we can choose $r$ large enough so that all sectional curvatures of 
$\mu_{s}'$ are are within $\epsilon $ of -1.
Note that the constructions of the Riemannian metrics
$\mu$ and $\mu_{0}'$ above commute with the
the cover $p:\bM\ra\M$, that is, $p^{*}\mu =\mu_{0}'$.

Since $f_{1/2}$ is the identity, $\M_{1/2}=\bar{M}$ and
$(\lambda_{r})_{1/2}$ is given by:
\begin{equation*}
(\lambda_{r})_{1/2}(x,u,t) = \begin{cases}
\bar{\rho}_{r}(x,u,t) &\text{$t\leq 2$, $5\leq t$}\\ 
& \\
\cosh^{2}(\alpha t)\sigma_{N}(x) + (\frac{e^{\alpha t}+
\delta_{1}(t)e^{-\alpha t}}{2})^{2} \sigma_{\bS^{1}}(u) + \alpha^2 dt^2 &\text{$2\leq t\leq
3$}\\
& \\
\cosh^{2}(\alpha t)\sigma_{N}(x) + \cosh^{2}(\alpha t) \sigma_{\bS^{1}}(u)+\alpha^2 dt^2
&\text{$3\leq t\leq 4$}\\
& \\
\cosh^{2}(\alpha t) \sigma_{N}(x) + (\frac{e^{\alpha t}+ \delta_{3}(t)e^{-\alpha t}}{2})^{2} 
\sigma_{\bS^{1}}(u)+ \alpha^2 dt^2 &\text{$4\leq t\leq 5$}\end{cases}
\end{equation*}
where $\alpha =r/6$.

We now define Riemannian metrics $\mu_{s}'$, $s\in [1/2,1]$.
Let $\eta : [1/2,1]\ra [1,2]$ be a smooth function such that
$\eta (1/2)=1$, $\eta (1)=0$ and $\eta$ is constant near
1/2 and 1. For $s\in [1/2,1]$, define
\begin{equation*}
(\lambda_{r})_{s}(x,u,t)=\begin{cases}
\bar{\rho}_{r}(x,u,t) &\text{$t\leq 2$, $5\leq t$}\\
\cosh^{2}(\alpha t) \sigma_{N}(x) + \Biggl(\frac{e^{\alpha t}+
\Bigl[\eta (s) (1+\delta_{1}(t))-1 \Bigr]  e^{-\alpha 
t}}{2} \Biggr)^{2} \sigma_{\bS^{1}}(u) + \alpha^2 dt^2 &\text{$2\leq t\leq 3$}\\
\cosh^{2}(\alpha t)\sigma_{N}(x) + \Biggl(\frac{e^{\alpha t}+
\Bigl[ 2\eta (s) - 1 \Bigr]e^{-\alpha t}}{2}\Biggr)^{2} 
\sigma_{\bS^{1}}(u) + \alpha^2 dt^2 &\text{$3\leq t\leq 4$}\\
\cosh^{2}(\alpha t)\sigma_{N}(x) + \Biggl( \frac{e^{\alpha t}+
\Bigl[\eta (s) (1+\delta_{3}(t))-1 \Bigr]e^{-\alpha 
t}}{2}\Biggr)^{2} \sigma_{\bS^{1}}(u) + \alpha^2 dt^2 &\text{$4\leq t\leq 5$}
\end{cases}
\end{equation*}
where $\alpha =r/6$.

Construct Riemannian metrics $\mu_{s}'$, $s\in [1/2,1]$ in the same way 
as before: pull back the metrics $(\lambda_{r})_{s}$ and fit 
them in $\bar{M}$.
Note that $(\lambda_{r})_{1}=\bar{\rho}_{r}$; hence $\mu_{1}'$ is the
hyperbolic metric on $\bar{M}$. Also, as before, we can suppose,
taking $r$ sufficiently large, that all sectional curvatures of 
$\mu_{s}'$ are are within $\epsilon $ of -1.

Finally, the smooth isotopy $f_{s}$, 
$s\in [0,1/2]$, induces a (top)  isotopy
$g_{t}:\bM\ra\bM$, with $g_{0}=1_{\bM}$ and $g_{s}:\bM\ra\M_{s}$
a diffeomorphism. Define $\mu_{s}=g_{s}^{*}\mu_{s}'$, for
$s\in [0,1/2]$, and $\mu_{s}=\mu_{s}'$, for $s\in [1/2,1]$.
It is straightforward to verify that we can choose the isotopy
$g_{s}$ in such a way that $\mu_{s}$ satisfies (iv) of the Proposition.

\section*{Appendix}

Here we sketch how to deduce the formula for the sectional curvature
of a doubly warped metric used to prove Lemma \ref{L:lem22}.

We have three steps. First we calculate the Levi-Civita connection. Then 
the curvature operator, and finally the sectional curvatures.

As in section 2, let $M_{1}$ and $M_{2}$
be Riemannian manifolds with Riemannian metrics $\sigma_{1}$ and 
$\sigma_{2}$, respectively. Let  also $\phi_{i} :\R\ra (0,\infty )$ 
be  smooth functions, $i=1,2.$ and define the doubly warped 
metric $\rho$ on $M=M_{1}\x M_{2}\x \R$:
\begin{equation*}
\rho = \phi_{1}^{2}\sigma_{1}+ \phi_{2}^{2}\sigma_{2}+ dt^{2}
\end{equation*}
or, equivalently,
\begin{equation*}
\langle u_{1}+v_{1}+s_{1} \p ,  u_{2}+v_{2}+s_{2} \p \rangle = 
\phi_{1}^{2}\langle u_{1},v_{1}\rangle_{1} + \langle 
u_{2},v_{2}\rangle_{2} + s_{1}s_{2}
\end{equation*}
where $\langle\ ,\ \rangle =\langle\ ,\ \rangle_{\rho}$, 
$\langle\ ,\ \rangle_{i}=\langle\ , \ \rangle_{\sigma_{i}}$,
$\p=\frac{\p}{\p t}$, $u_{i}+v_{i}+s_{i}\p \in T_{x_{1}}M_{1}
\oplus T_{x_{2}}M_{2}\oplus\R =T_{(x_{1},x_{2},t)}M$.

Note that if  $u_{1},u_{2}$ are vector fields
on $M$ which are zero in the $M_{2}$ and $\R$ directions,
then $\langle u_{1},u_{2}\rangle = \phi_{1}^{2}\langle u_{1},u_{2}\rangle_{1}$.
Analogously, for vector fields $v_{1}$, $v_{2}$
on $M$ which are zero in the $M_{1}$ and $\R$ directions, we have
$\langle v_{1},v_{2}\rangle = \phi_{2}^{2}\langle v_{1},v_{2}\rangle_{2}$.

{\bf The Connection.}

We use the Koszul formula that relates the Levi-Civita connection $D$ with the 
Riemannian metric $(\ ,\ )$ of a Riemannian manifold:
\begin{equation}
2(Z, D_{Y}X) = X(Y,Z)+Y(X,Z)-Z(X,Y)    
\label{E:(1)}
\end{equation}
Here we are assuming that the vector fields $X,Y,Z$ commute. 

Now $u, u_{1},u_{2},\ldots$ will denote tangent 
vectors in $TM_{1}\subset TM$, or vector fields
on $M$ which are constant in the $M_{2}$ and $\R$ directions.
Analogously,  $v, v_{1},v_{2},\ldots$ will denote tangent 
vectors in $TM_{2}\subset TM$, or vector fields
on $M$ which are constant in the $M_{1}$ and $\R$ directions.
We assume that all vector fields commute.
Let $\n, \n^{1},\n^{2}$ denote the 
Levi-Civita connections of the Riemannian manifolds $M$, $M_{1}$,
$M_{2}$, respectively. Write $\p =\frac{\p}{\p t}$.

\begin{Cla}\label{C:cla1}
\begin{equation}
\begin{aligned}
(i) &\quad \n_{\p} \p = 0\\
(ii) &\quad \n_{\p}  u = \n_{u} \p  = \frac{\phi_{1}'}{\phi_{1}} u\\
(iii) &\quad \n_{\p} v  = \n_{v} \p  = \frac{\phi_{2}'}{\phi_{2}} v\\
(iv) &\quad \n_{u_{1}} u_{2} = -\frac{\phi_{1}'}{\phi_{1}} \langle 
u_{1},u_{2}\rangle \p +\n^{1}_{u_{1}}u_{2}\\
(v) &\quad \n_{v_{1}} v_{2} = -\frac{\phi_{2}'}{\phi_{2}}
\langle v_{1},v_{2}\rangle \p + \n^{2}_{v_{1}}v_{2}\\
(vi) &\quad \n_{u}v  = \n_{v}u  = 0
\end{aligned}
\label{E:(2)}
\end{equation}
\end{Cla}

\begin{proof}
(i) By (\ref{E:(1)}) above, we have $2\langle u , \n_{\p}\p\rangle=
2\langle v ,\n_{\p}\p \rangle=2\langle \p , \n_{\p}\p\rangle=0$.
Hence $\n_{\p}\p =0$.

(ii) By (\ref{E:(1)}) we have $ \langle v , \n_{u}\p\rangle =0$.
Since $\langle \p , \p\rangle =1$, also by (\ref{E:(1)}), we have 
$\langle \p , \n_{u}\p\rangle =0$.

Finally, again by (\ref{E:(1)}), 
$2\langle u_{1} , \n_{u}\p\rangle =\p \langle u_{1} , u\rangle$.

Hence $\langle u_{1} , \n_{u}\p\rangle=\frac{1}{2}\p\langle u_{1} , 
u\rangle= \frac{1}{2}\p[\phi_{1}^{2}
\langle u_{1} , u\rangle_{1}]=\phi_{1}'\phi_{1}
\langle u_{1} , u\rangle_{1} = \frac{\phi_{1}'}{\phi_{1}}\langle u_{1} 
,u \rangle 
=\langle u_{1} ,\frac{\phi_{1}'}{\phi_{1}} u\rangle$.

It follows that $\n_{u}\p=\frac{\phi_{1}'}{\phi_{1}} u$.

(iii) Same as (ii).

(iv) Since $v\langle u_{1} , u_{2}\rangle =0$, we have, by (\ref{E:(1)}),
that $\langle v, \n_{u_{1}}u_{2} \rangle =0$, for all $v$.

Also, as in the proof of (ii), we have $\langle \p , \n_{u_{1}}u_{2} \rangle  
=-\frac{\phi_{1}'}{\phi_{1}}\langle u_{1} , u_{2} \rangle$.

Finally, for all $u_{3}$ we have
$2\langle u_{3} , \n_{u_{1}}u_{2} \rangle =u_{2} \langle u_{1} , u_{3} \rangle 
+u_{1} \langle u_{2} , u_{3} \rangle- u_{3} \langle u_{1} , u_{2} \rangle=
u_{2} \phi_{1}^{2}\langle u_{1} , u_{3} \rangle_{1} 
+u_{1}\phi_{1}^{2} \langle u_{2} , u_{3} \rangle_{1}- u_{3}
\phi_{1}^{2}\langle u_{1} , u_{2} \rangle_{1}=\phi_{1}^{2}\Bigl\{
u_{2} \langle u_{1} , u_{3} \rangle_{1} 
+u_{1} \langle u_{2} , u_{3} \rangle_{1}- u_{3} \langle u_{1} , u_{2} 
\rangle_{1}\Bigr\} =\phi_{1}^{2}\Bigl\{ 2\langle u_{3} , 
\n^{1}_{u_{1}}u_{2} \rangle_{1}\Bigr\} =
2\langle u_{3} , \n^{1}_{u_{1}}u_{2} \rangle$.
This proves (iv).

(v) Same as (iv).

(vi) Since $u\langle v , v_{1} \rangle =  
v\langle u , u_{1} \rangle =0$,
by (\ref{E:(1)}) we have $\langle u_1 , \n_{u}v \rangle = \langle v_1 , 
\n_{u}v \rangle = \langle \p , \n_{u}v \rangle=0$.  It follows that 
$\n_{u}v=0$. This proves the claim.  
\end{proof}

{\bf The Curvature Operator.}

Let $R$ denote the curvature tensor (of type (3,1)) of the Riemannian 
manifold $M$; that is $R_{ab}c = \n_b\n_a - \n_a\n_b-\n_{[b,a]}c$.
A straightforward calculation, using (\ref{E:(2)}), shows:
\begin{equation}
\begin{aligned}
R_{\p u} \p &= -\frac{\phi_{1}''}{\phi_{1}} u\\
R_{\p v} \p &= -\frac{\phi_{2}''}{\phi_{2}} v\\
R_{u_{1}u_{2}} \p &= 0\\    
R_{v_{1}v_{2}} \p &= 0\\  
R_{uv}\p &= 0
\end{aligned}
\label{E:(3)}
\end{equation}    
and
\begin{equation}\begin{aligned}
R_{v\p} u &= 0 &&\quad  R_{u \p} v = 0\\    
R_{v_{1}v_{2}} u &= 0 &&\quad R_{u_{1}u_{2}} v = 0\\      
R_{u_{1}\p} u_{2} &= -\frac{\phi_{1}''}{\phi_{1}} \langle 
u_{1},u_{2}\rangle \p  &&\quad R_{v_{1}\p} v_{2} = -\frac{\phi_{2}''}{\phi_{2}} \langle 
v_{1},v_{2}\rangle  \p\\
R_{u_{1}v} u_{2} &= -\frac{\phi_{1}'\phi_{2}'}{\phi_{1}\phi_{2}} \langle u_{1} , 
u_{2}\rangle  v &&\quad R_{v_{1}u} v_{2} =-\frac{\phi_{1}'\phi_{2}'}{\phi_{1}\phi_{2}}
 \langle v_{1} , v_{2}\rangle u
\end{aligned}
\label{E:(4)}
\end{equation}    
Note that (\ref{E:(1)}) implies:
\begin{equation*}
\begin{cases}
u_{2}\langle u_{1}, u_{3}\rangle+\langle u_{2}, \n_{u_{1}}u_{3}\rangle=  
u_{1}\langle u_{2}, u_{3}\rangle+\langle u_{1}, 
\n_{u_{2}}u_{3}\rangle\\
v_{2}\langle v_{1}, v_{3}\rangle+\langle v_{2}, \n_{v_{1}}v_{3}\rangle=  
v_{1}\langle v_{2}, v_{3}\rangle +\langle v_{1}, 
\n_{v_{2}}v_{3}\rangle\end{cases}
\end{equation*}
>From this and (\ref{E:(2)}) it follows that
\begin{equation}
\begin{aligned}
R_{u_{1}u_{2}}u_{3} &=R^{1}_{u_{1}u_{2}} u_{3} -
\frac{(\phi_{1}')^{2}}{\phi_{1}^{2}}\Bigl(\langle u_{1} , u_{3}\rangle
u_{2} - \langle u_{2},u_{3} \rangle u_{1}\Bigr)\\
R_{v_{1}v_{2}} v_{3} &= R^{2}_{v_{1}v_{2}} v_{3} -
\frac{(\phi_{2}')^{2}}{\phi_{2}^{2}}\Bigl(\langle v_{1} , v_{3}\rangle
v_{2} - \langle v_{2},v_{3} \rangle v_{1}\Bigr)
\end{aligned}
\label{E:(5)}
\end{equation}

Let $\mR$ denote the curvature operator on two-forms for $M$, defined by
$\langle \mR (a\wedge b), c\wedge d \rangle = \langle R_{ab}c, 
d\rangle$. 
Also, let $\mR^{i}$ denote the curvature operator on two-forms for 
$M_{i}$.
Recall that the scalar product on $\wedge^{2}TM$ is given 
by $\langle a\wedge b,c\wedge d \rangle =\langle a, c \rangle \langle 
b, d\rangle - 
\langle a, d\rangle \langle b,c\rangle.$

\begin{Cla}\label{C:cla2}
\begin{equation}
\begin{aligned}
(i) &\quad \mR (\partial \wedge u) = -\frac{\phi_{1}''}{\phi_{1}}\p\wedge u\\
(ii) &\quad \mR (\partial \wedge v) = -\frac{\phi_{2}''}{\phi_{2}}\p\wedge v\\
(iii) &\quad \mR (u_{1}\wedge u_{2}) = \mR^{1} (u_{1}\wedge u_{2}) - 
\frac{(\phi_{1}')^{2}}{\phi_{1}^{2}}  u_{1}\wedge u_{2}\\
(iv) &\quad \mR (v_{1}\wedge v_{2}) = \mR^{2} (v_{1}\wedge v_{2}) - 
\frac{(\phi_{2}')^{2}}{\phi_{2}^{2}}  v_{1}\wedge v_{2}\\
(v) &\quad \mR (u\wedge v) =  -\frac{\phi_{1}'\phi_{2}'}{\phi_{1}\phi_{2}}   u\wedge v
\end{aligned}
\label{E:(6)}
\end{equation}
\end{Cla}

\begin{proof} (i) From (\ref{E:(3)}) or (\ref{E:(4)}) we  have $\langle \mR (\p\wedge u), 
\p\wedge v\rangle =0$. From (\ref{E:(4)}) we have $\langle \mR (\p\wedge u), 
v_{1}\wedge v_{2}\rangle = \langle R(\partial \wedge u),u_1\wedge u_2\rangle = \langle
R(\partial\wedge u,u_1\wedge v)\rangle =0$.
Hence $\mR(\p\wedge u)$ is a linear 
combination of two-vectors of the form $\p\wedge u_{1}$.
>From (\ref{E:(3)}) or (\ref{E:(4)}) we  have that, for all
$u_{1}$,  $\langle \mR (\p\wedge u), 
\p\wedge u_{1}\rangle =-\frac{\phi_{1}'}{\phi_{1}}\langle 
u,u_{1}\rangle =\langle -\frac{\phi_{1}'}{\phi_{1}}\p\wedge u, 
\p\wedge u_{1}\rangle $. This proves (i).

(ii) Same as (i).

(iii) By (\ref{E:(3)}),  $\langle \mR (u_{1}\wedge u_{2}), 
\p\wedge v\rangle = \langle R(u_1\wedge u_2),\partial\wedge u\rangle =0$.
By (\ref{E:(4)}),  $\langle \mR (u_{1}\wedge u_{2}), 
v_{1}\wedge v_{2}\rangle = \langle R(u_1\wedge u_2),u\wedge v\rangle=0$.
Hence $\mR(u_{1}\wedge u_{2})$ is a linear 
combination of two-vectors of the form $u_{i}\wedge u_{j}$.
But, by (\ref{E:(5)}), we have 
\begin{equation*}
\begin{aligned} 
\langle \mR (u_{1}\wedge u_{2}), u_{3}\wedge u_{4}\rangle
&=\langle \mR^{1} (u_{1}\wedge u_{2}), u_{3}\wedge u_{4}\rangle
-\frac{(\phi_{1}')^{2}}{\phi_{1}^{2}}
\Bigl(  \langle u_{1},u_{3}\rangle \langle u_{2},u_{4}\rangle-
\langle u_{2},u_{3}\rangle\langle u_{1},u_{4}\rangle
\Bigr)\\
&=\langle \mR^{1} (u_{1}\wedge u_{2}), u_{3}\wedge u_{4}\rangle
-\frac{(\phi_{1}')^{2}}{\phi_{1}^{2}}
\langle u_{1}\wedge u_{2},u_{3}\wedge u_{4}\rangle\\
&=\langle \Bigl\{ \mR^{1} (u_{1}\wedge u_{2}) -
\frac{(\phi_{1}')^{2}}{\phi_{1}^{2}} u_{1}\wedge u_{2}\Bigl\}
, u_{3}\wedge u_{4}\rangle
\end{aligned}
\end{equation*}
for all $u_{3}$, $u_{4}$. This proves (iii).

(iv) Same as (iii).

(v) By (\ref{E:(3)}) and (\ref{E:(4)}), $\mR (u\wedge v)$ is a linear combination
of two-vectors of the form $u_{i}\wedge v_{j}$. From (\ref{E:(4)}) we 
have that
\begin{equation*}
\langle \mR (u\wedge v), u_{1}\wedge v_{1}\rangle
=- \frac{\phi_{1}'\phi_{2}'}{\phi_{1}\phi_{2}}
\langle u , u_{1}\rangle  \langle v , v_{1}\rangle
=\langle  - \frac{\phi_{1}'\phi_{2}'}{\phi_{1}\phi_{2}}
 u\wedge v , u_{1}\wedge v_{1}\rangle
\end{equation*}
for all $u_{1}$, $v_{1}$. This proves (v) and the claim.
\end{proof}

{\bf The Sectional Curvature.}

Let $P\subset T_{(x_{1}, x_{2},t)}(M_{1}\x M_{2}\x\R ) = 
T_{x_{1}}M_{1}\oplus T_{x_{2}}M_{2}\oplus \R $ be the two-plane 
generated by the orthonormal basis $\{ u_{1}+v_{1}+s\frac{\partial}{\partial 
t} , u_{2}+v_{2}\}$, where $u_{1},u_{2}\in T_{x_{1}}M_{1}$,
$v_{1},v_{2}\in T_{x_{2}}M_{2}$. 

Write $a=u_{1}+v_{1}+s \p$ and $b=u_{2}+v_{2}$.
Then the sectional curvature $K(P)$ is given by $\langle \mR (a,b),a\wedge b\rangle$.
Note that 
\begin{equation*}
\begin{cases}
\mR(a\wedge b)=\mR(u_{1}\wedge u_{2})+ \mR(u_{1}\wedge v_{2}) +
\mR(v_{1}\wedge u_{2})+\mR(
v_{1}\wedge v_{2})+ s\mR(\p\wedge u_{2})+ s\mR(\p\wedge v_{2})\\
a\wedge b= u_{1}\wedge u_{2}+ u_{1}\wedge v_{2} +v_{1}\wedge u_{2}+
v_{1}\wedge v_{2}+ s\p\wedge u_{2}+ s\p\wedge v_{2}
\end{cases}.
\end{equation*}
Multiplying these terms, a straightforward
calculation using (\ref{E:(6)}) shows that our formula for the
sectional curvature of a doubly warped metric holds.

\sk1

F.T. Farrell

SUNY, Binghamton, N.Y., 13902, U.S.A.
\sk1

P. Ontaneda

UFPE, Recife, PE 50670-901, Brazil


\begin{thebibliography}{99}

\bibitem{A} S.I. Al'ber, {\em Spaces of mappings into manifold of negative
curvature}, Dokl. Akad. Nauk. SSSR {\bf 178} (1968), 13-16.
    
\bibitem{BO} R.L. Bishop and B. O'Neill, {\em Manifolds of negative 
curvature}, Trans. Amer. Math. Soc. {\bf145} (1969) 1-49. 
   
\bibitem{BL} D. Burghelea and R. Lashof, {\em Stability of 
concordances and suspension homeomorphism}, Ann. of Math. (2) {\bf 
105} (1977), 449-472.

\bibitem{BK} K. Burns and A. Katok, {\em Manifolds with non-positive curvature}, Ergodic
Theory \& Dynam. Sys. {\bf 5} (1985), 307-317.

\bibitem{Ce} J. Cerf, {\em La stratification naturelle des espaces de fonctions differentiables
reels et le theoreme de la pseudo-isotopie}, Inst. Hautes Etudes Sci. Publ. Math. No. {\bf 39}
(1970), 5-173.

\bibitem{Co} K. Corlette, {\em Archimedean superrigidity and harmonic geometry}, Ann. of Math.
{\bf 135} (1992) 165-182.

\bibitem{E} R.D. Edwards, {\em The topology of manifolds and cell-like maps}, in Proc. of the
ICM (Helsinki, 1978), pp. 111-127, Acad. Sci. Fennica, Helsinki, 1980.

\bibitem{EK} R.D. Edwards and R.C. Kirby, {\em Deformations of spaces of imbeddings}, Ann. of
Math. {\bf 93} (1971), 63-88.

\bibitem{EL1} J. Eells and L. Lemaire, {\em A report on harmonic maps}, Bull. of LMS, {\bf 10}
(1978), 1-68.

\bibitem{EL2} J. Eells and L. Lemaire, {\em Deformations of metrics and associated harmonic
maps}, Patodi Memorial Vol., Geometry and Analysis (Tata Inst., 1981) 33-45.

\bibitem{EL3} J. Eells and L. Lemaire, {\em Selected topics in harmonic maps}, CBMS
Regional Conf. Series {\bf 50}, Amer. Math. Soc., Providence, R.I, 1983.

\bibitem{EL4} J. Eells and L. Lemaire, {\em Another report on harmonic maps}, Bull. of LMS,
{\bf 20} (1988), 385-524.

\bibitem{ES} J. Eells and J.H. Sampson, {\em Harmonic mappings of Riemannian
manifolds}, Amer. J. Math. {\bf 86} (1964), 109-160.

\bibitem{FJ1} F.T. Farrell and L.E. Jones, {\em K-theory and dynamics II},
Ann. of Math. {\bf 126} (1987) 451-493.

\bibitem{FJ2} F.T. Farrell and L.E. Jones, {\em Negatively curved manifolds with exotic smooth
structures}, J. Amer. Math. Soc. {\bf 2} (1989) 899-908.

\bibitem{FJO} F.T. Farrell, L.E. Jones and P. Ontaneda, {\em Examples of non-homeomorphic
harmonic maps between negatively curved manifolds}, Bull. London Math. Soc. {\bf 30} 
(1998) 295-296.

\bibitem{FOR} F.T. Farrell, P. Ontaneda and M.S. Raghunathan, 
{\em Non-univalent harmonic maps homotopic to diffeomorphisms},
Jour. Diff. Geom. {\bf 54} (2000)    227-253.

\bibitem{GS} M. Gromov and R. Schoen, {\em Harmonic maps into singular spaces and $p$-adic
superrigidity of lattices in groups of rank one}, Inst. Hautes {\'E}tudes Sci. Publ. Math. 
{\bf 76} (1992) 165-246.

\bibitem{Har} P. Hartman, {\em On homotopic harmonic maps}, Canad. J. 
Math. {\bf19}, (1967) 673-687.

\bibitem{Hat} A.E. Hatcher, {\em Concordance spaces, higher simple homotopy theory, and
applications}, Proc. Symp. Pure Math. {\bf 32} (1978), 3-21.

\bibitem{I} K. Igusa, {\em Stability theorems for pseudoisotopies}, 
K-theory {\bf 2} (1988), 1-355.

\bibitem{JY} J. Jost and S.-T. Yau, {\em Harmonic maps and superrigidity}, Proc. Sympos. Pure
Math {\bf 54} (Amer. Math. Soc., Providence, R.I., 1993), 245-280.

\bibitem{MSY} N. Mok, Y.-T. Siu and S.-K. Yeung, {\em Geometric superrigidity}, Invent. Math.
{\bf 113} (1993) 57-83.

\bibitem{M} G.D. Mostow, {\em Quasi-conformal mappings in $n$-space    
and the rigidity of hyperbolic space forms}, Inst. Hautes \'Etudes
Sci. Publ. Math. {\bf 34} (1967), 53-104.
    
\bibitem{O} P. Ontaneda, {\em Hyperbolic manifolds with negatively curved 
exotic triangulations in dimension six},  J. Diff. Geom. {\bf 40} (1994), 7-22. 

\bibitem{Sa} J. Sampson, {\em Some properties and applications of harmonic mappings}, Ann. Sci.
{\'E}cole Norm. Sup. {\bf 11} (1978) 211-228.

\bibitem{Sc} M. Scharlemann, {\em Smooth CE maps and smooth homeomorphisms}, Lecture Notes in
Math., vol. 341, Springer-Verlag, Berlin and New York, 1978, pp. 234-240.

\bibitem{SS} M. Scharlemann and L. Siebenmann, {\em The 
Hauptvermutung for smooth singular homeomorphisms}, in Manifolds 
Tokyo 1973, Akio Hattori ed.
Univ. of Tokyo Press, 85-91.

\bibitem{SY} R. Schoen and S.-T. Yau, {\em On univalent harmonic maps between surfaces},
Invent. Math. {\bf 44} (1978), 265-278.

\bibitem{Si} L.C. Siebenmann, {\em Approximating cellular maps by homeomorphisms}, Topology 
{\bf 11} (1973), 271-294.

\bibitem{Sui} Y.-T. Siu, {\em The complex-analyticity of harmonic maps and the strong rigidity
of compact K{\"a}hler manifolds}, Ann. of Math. {\bf 112} (1980) 73-111.

\end{thebibliography}
\end{document}